# On Fredholm's Integral Equations on the Real Line, Whose Kernels Are Linear in a Parameter

Igor M. Novitskii
Khabarovsk Division
Institute of Applied Mathematics
Far-Eastern Branch of the Russian Academy of Sciences
54, Dzerzhinskiy Street, Khabarovsk 680 000, RUSSIA
Email: novim@iam.khv.ru

***Abstract*—In this paper, we study an infinite system of Fredholm series of polynomials in $\lambda$, formed, in the classical way, for a continuous Hilbert-Schmidt kernel on $\mathbb{R}\times\mathbb{R}$ of the form $\boldsymbol{H}(s,t)-\lambda\boldsymbol{S}(s,t)$, where $\lambda$ is a complex parameter. We prove a convergence of these series in the complex plane with respect to sup-norms of various spaces of continuous functions vanishing at infinity. The convergence results enable us to solve explicitly an integral equation of the second kind in $L^2(\mathbb{R})$, whose kernel is of the above form, by mimicking the classical Fredholm-determinant method.**

## I. Introduction and Preliminaries

The higher Fredholm minors play an essential part in Fredholm's [2] solution of the integral equations of the second kind in the space $C[0,1]$, whose kernels are defined and continuous on $[0,1]^2$; the objects provide us with explicit bases for the null spaces and explicit general solutions, for all values of the spectral parameter. The Carleman-Mikhlin-Smithies theory (see [1], [9], [16]) dealing with any (measurable) Hilbert-Schmidt kernel of possibly unbounded support extends to $L^2$ integral operators only those parts of the Fredholm theory which concern the Fredholm determinant and the first Fredholm minor. Among other things in the present paper we have managed to construct all of the $p$-th Fredholm minors ($p=1,2,3,\dots$) for Hilbert-Schmidt kernels on $\mathbb{R}^2$, but at the price of imposing on the kernels the extra condition of being in the class of the so called $K^0$ *kernels of Mercer type* (see Subsections I-3 and II-5). Luckily, this condition can always be achieved by means of a unitary equivalence transformation (see Theorem 2 below); and this is, therefore, not a serious loss of generality when working in that class of kernels.

Before we can write down and prove our results, we need to fix the terminology and notation and to give some definitions and foregoing results.

*1) Spaces:* Throughout, $\mathcal{H}$ is a complex, separable, infinite-dimensional Hilbert space with norm $\|\cdot\|_{\mathcal{H}}$ and inner product $\langle\cdot,\cdot\rangle_{\mathcal{H}}$. $\mathbb{R}$ is the real line equipped with the Lebesgue measure, and $L^2=L^2(\mathbb{R})$ is the Hilbert space of (equivalence classes of) measurable complex-valued functions on $\mathbb{R}$ equipped with the inner product $\langle f,g\rangle=\int f(s)\overline{g(s)}\,ds$ and the norm $\|f\|=\langle f,f\rangle^{\frac{1}{2}}$. (Integrals with no indicated domain, such as above, are always to be extended over $\mathbb{R}$.) $C(X,B)$, where $B$ is a Banach space with norm $\|\cdot\|_B$, is the Banach algebra (with the norm $\|f\|_{C(X,B)}=\sup_{x\in X}\|f(x)\|_B$) of continuous $B$-valued functions defined on a locally compact space $X$ and *vanishing at infinity* (that is, given any $f\in C(X,B)$ and $\varepsilon>0$, there exists a compact subset $X(\varepsilon,f)\subset X$ such that $\|f(x)\|_B<\varepsilon$ whenever $x\notin X(\varepsilon,f)$). A series $\sum_n f_n$ is $B$-*absolutely convergent in* $C(X,B)$ if $f_n\in C(X,B)$ for all $n$ and the series $\sum_n\|f_n(x)\|_B$ converges in $C(X,\mathbb{R})$. $\theta_p$ denotes the zero element in the space $C\left(\mathbb{R}^{2p},\mathbb{C}\right)$. $\mathbb{C}$ is the complex plane. When $n=0$, $C\left(\mathbb{R}^n,B\right)$ is identified with $B$.

*2) Linear Operators:* Throughout, $\mathfrak{R}(\mathcal{H})$ is the Banach algebra of all bounded linear operators acting on $\mathcal{H}$. If $T\in\mathfrak{R}(\mathcal{H})$, then $T^*$ stands for the adjoint to $T$ (w.r.t. $\langle\cdot,\cdot\rangle_{\mathcal{H}}$). An operator $P\in\mathfrak{R}\left(L^2\right)$ is *positive* (written $P\geqslant 0$) if $\langle Px,x\rangle\geqslant 0$ for all $x\in L^2$. If $A=A^*$, $B=B^*\in\mathfrak{R}\left(L^2\right)$, then we write $A\leqslant B$ iff $B-A\geqslant 0$. The operator family $\mathcal{M}^+(T)$ where $T\in\mathfrak{R}\left(L^2\right)$ is the set of all those positive operators $P\in\mathfrak{R}\left(L^2\right)$ that are factorizable as $P=TB$ or as $P=BT$ with $B\in\mathfrak{R}\left(L^2\right)$. A factorization of $T\in\mathfrak{R}\left(L^2\right)$ into the product $T=WV^*$ ($V$, $W\in\mathfrak{R}\left(L^2\right)$) is called an $\mathcal{M}$ *factorization* for $T$ provided that $VV^*$, $WW^*\in\mathcal{M}^+(T)$.

*Example 1:* One example of an $\mathcal{M}$ factorization for any $T\in\mathfrak{R}\left(L^2\right)$ is obtained by letting $W=UP$, $V=P$, where $P$ is the positive square root of $|T|:=(T^*T)^{\frac{1}{2}}$ and $U$ is the partially isometric factor in the polar decomposition $T=U|T|$. Indeed: $T=WV^*$, $WW^*=U|T|U^*=TU^*=UT^*\in\mathcal{M}^+(T)$, and $VV^*=|T|=T^*U=U^*T\in\mathcal{M}^+(T)$.

A bounded linear operator $U\colon\mathcal{H}\to L^2$ is *unitary* if it has range $L^2$ and $\langle Uf,Ug\rangle=\langle f,g\rangle_{\mathcal{H}}$ for all $f$, $g\in\mathcal{H}$.

*3) Integral operators:* An operator $T\in\mathfrak{R}\left(L^2\right)$ is *integral* if there is a complex-valued measurable function $\boldsymbol{T}$ (*kernel*) defined on the Cartesian product $\mathbb{R}^2=\mathbb{R}\times\mathbb{R}$ such that, for each $f$ in $L^2$, $(Tf)(s)=\int\boldsymbol{T}(s,t)f(t)\,dt$ for a.e. $s$ in $\mathbb{R}$ [4], [5]. A function $\boldsymbol{T}\in C\left(\mathbb{R}^2,\mathbb{C}\right)$ is called a $K^0$ kernel if the so-called *Carleman functions* $\boldsymbol{t}$, $\boldsymbol{t'}\colon\mathbb{R}\to L^2$, defined via $\boldsymbol{T}$ by $\boldsymbol{t}(s)=\overline{\boldsymbol{T}(s,\cdot)}$, $\boldsymbol{t'}(s)=\boldsymbol{T}(\cdot,s)$, belong to $C\left(\mathbb{R},L^2\right)$. A $K^0$ kernel $\boldsymbol{T}$ is called of *Mercer type* if it induces an integral operator $T\in\mathfrak{R}\left(L^2\right)$, with the property that any operator belonging to $\mathcal{M}^+(T)$ is also an integral operator induced by a $K^0$ kernel. A counter-example in [10, Section 2] shows that there might be $K^0$ kernels which are not of Mercer type.

The following theorem can be seen as a generalization of Mercer's [8] theorem (about absoluteness and uniformity of convergence of bilinear eigenfunction expansions for continuous compactly supported kernels of positive, integral operators) to various other settings (for a detailed discussion see [10], [13]):

*Theorem 1:* Let $T \in \mathfrak{R}\left(L^2\right)$ be an integral operator, with a kernel $\boldsymbol{T}$ that is a $K^0$ kernel of Mercer type. Then, for any $\mathcal{M}$ factorization $T = WV^*$ for $T$ and for any orthonormal basis $\{u_n\}$ for $L^2$, the following bilinear formula holds

$$\boldsymbol{T}(s,t) = \sum_n Wu_n(s)\overline{Vu_n(t)} \quad \text{for all } s,\, t \in \mathbb{R},$$

where the series converges $\mathbb{C}$-absolutely in $C\left(\mathbb{R}^2, \mathbb{C}\right)$.

The next theorem says that integral operators that have as a kernel a $K^0$ kernel of Mercer type are not so rare as they might appear to be.

*Theorem 2:* Suppose that for an operator family $\{S_\gamma\}_{\gamma\in\mathcal{G}} \subset \mathfrak{R}(\mathcal{H})$ with an index set of arbitrary cardinality there exists an orthonormal sequence $\{e_n\}$ in $\mathcal{H}$ such that

$$\lim_{n\to\infty} \sup_{\gamma\in\mathcal{G}} \|S_\gamma e_n\|_{\mathcal{H}} = 0, \quad \lim_{n\to\infty} \sup_{\gamma\in\mathcal{G}} \left\|(S_\gamma)^* e_n\right\|_{\mathcal{H}} = 0.$$

Then there exists a unitary operator $U : \mathcal{H} \to L^2$ such that all the operators $T_\gamma = US_\gamma U^{-1}$ $(\gamma \in \mathcal{G})$ and their linear combinations are integral operators on $L^2$, whose kernels are $K^0$ kernels of Mercer type.

For proofs and applications of this theorem see [10], [13], and [14].

An operator $A \in \mathfrak{R}\left(L^2\right)$ is *nuclear* if $\sum_n |\langle Au_n, u_n\rangle| < +\infty$ for any choice of an orthonormal basis $\{u_n\}$ of $L^2$. A *Hilbert-Schmidt kernel* $\boldsymbol{T}$ on $\mathbb{R}^2$ is a kernel for which $\int\int |\boldsymbol{T}(s,t)|^2\, dt\, ds < +\infty$. A nuclear operator is always an integral operator, whose kernel is Hilbert-Schmidt but need not be continuous.

*4) A Second-kind Integral Equation:* Consider the integral equation of the second kind in $L^2$

$$f(s) - \lambda \int \boldsymbol{S}(s,t) f(t)\, dt = g(s) \quad \text{for a.e. } s \text{ in } \mathbb{R}, \qquad (1)$$

where $\boldsymbol{S}\colon \mathbb{R}^2 \to \mathbb{C}$ (the *kernel* of the equation) is a given $K^0$ kernel of Mercer type inducing a nuclear operator $S$ on $L^2$, the scalar $\lambda \in \mathbb{C}$ (a *parameter*) is given, the function $g$ of $L^2$ is given, and the function $f$ of $L^2$ is to be determined. In [11] it is shown that equation (1) can be dealt with by means of a technique completely similar to that by Fredholm (see [2], [6]) for the same kind of integral equations, but with compactly supported continuous kernels. A brief outline of the technique is as follows. For $p = 0, 1, 2, \ldots$, write the (power) *Fredholm series* for $\boldsymbol{S}$ as

$$\sum_{n=0}^{\infty} \frac{(-\lambda)^n}{n!} \int \cdots \int \boldsymbol{S}\begin{pmatrix} s_1 \ldots s_p\ \xi_1 \ldots \xi_n \\ t_1 \ldots t_p\ \xi_1 \ldots \xi_n \end{pmatrix} d\xi_1 \ldots d\xi_n, \qquad (2)$$

where $\boldsymbol{S}\begin{pmatrix} x_1 \ldots x_\nu \\ y_1 \ldots y_\nu \end{pmatrix} := \begin{cases} \det\left[\boldsymbol{S}(x_i, y_j)\right]_{i,j=1}^{\nu} & \text{if } \nu > 0, \\ 1 & \text{if } \nu = 0 \end{cases}$ are the so-called *$\nu$-th compounds* of $\boldsymbol{S}$. The sum of the series (2) is the $p$-th *Fredholm minor* (written $\boldsymbol{S}\left(\begin{smallmatrix} s_1 \ldots s_p \\ t_1 \ldots t_p \end{smallmatrix} \middle| \lambda\right)$) when $p > 0$ and the *Fredholm determinant* (written $\boldsymbol{S}(\lambda)$) when $p = 0$. For fixed $\lambda \in \mathbb{C}$, every $p$-th Fredholm minor, thought of as a function from $\mathbb{R}^{2p}$ to $\mathbb{C}$, is in $C\left(\mathbb{R}^{2p}, \mathbb{C}\right)$, because its series may be shown to be uniformly convergent on $\mathbb{R}^{2p}$. When the parameter value $\lambda$ is given, there is a least non-negative integer $\mathbf{d} = \mathbf{d}(\lambda)$ such that

$$\boldsymbol{S}\left(\begin{matrix} s_1 \ldots s_{\mathbf{d}} \\ t_1 \ldots t_{\mathbf{d}} \end{matrix} \middle| \lambda\right) \neq \theta_{\mathbf{d}}. \qquad (3)$$

If points $s_1, \ldots, s_{\mathbf{d}}, t_1, \ldots, t_{\mathbf{d}}$ of $\mathbb{R}$ are chosen so that the above minor does not vanish, then the functions $v_1, \ldots, v_{\mathbf{d}}$, given by the equations

$$v_i(s) = \boldsymbol{S}\left(\begin{matrix} s_1 \ldots s_{i-1}\ s\ s_{i+1} \ldots s_{\mathbf{d}} \\ t_1 \ldots t_{i-1}\ t_i\ t_{i+1} \ldots t_{\mathbf{d}} \end{matrix} \middle| \lambda\right) \quad (i = 1, \ldots, \mathbf{d})$$

for all $s$ in $\mathbb{R}$, form a basis for the set of solutions $f$ of the homogeneous equation $f(s) - \lambda \int \boldsymbol{S}(s,t) f(t)\, dt = 0$, and similarly the functions $u_1, \ldots, u_{\mathbf{d}}$, given by the equations

$$u_i(t) = \overline{\boldsymbol{S}\left(\begin{matrix} s_1 \ldots s_{i-1}\ s_i\ s_{i+1} \ldots s_{\mathbf{d}} \\ t_1 \ldots t_{i-1}\ t\ t_{i+1} \ldots t_{\mathbf{d}} \end{matrix} \middle| \lambda\right)} \quad (i = 1, \ldots, \mathbf{d})$$

for all $t$ in $\mathbb{R}$, form a basis for the set of solutions $f$ of the conjugate homogeneous equation $f(t) - \overline{\lambda} \int \overline{\boldsymbol{S}(s,t)} f(s)\, ds = 0$. The equation (1) is solvable if and only if $\langle u_i, g\rangle = 0$ for $i = 1, \ldots, \mathbf{d}$, and the general solution is then given by

$$f(s) = g(s) + \lambda \int \frac{\boldsymbol{S}\left(\begin{matrix} s\ s_1 \ldots s_{\mathbf{d}} \\ t\ t_1 \ldots t_{\mathbf{d}} \end{matrix} \middle| \lambda\right)}{\boldsymbol{S}\left(\begin{matrix} s_1 \ldots s_{\mathbf{d}} \\ t_1 \ldots t_{\mathbf{d}} \end{matrix} \middle| \lambda\right)} g(t)\, dt + \sum_{i=1}^{\mathbf{d}} c_i v_i(s)$$

where $c_1, \ldots, c_{\mathbf{d}}$ are arbitrary complex constants.

## II. Main Results

*5) Assumptions:* Equation (1) is a particular case of the integral equation

$$f(s) + \int \left(\boldsymbol{H}(s,t) - \lambda \boldsymbol{S}(s,t)\right) f(t)\, dt = g(s) \text{ a.e. in } \mathbb{R}, \quad (4)$$

in which $\boldsymbol{S}$ is as in (1), and the added kernel $\boldsymbol{H} \neq \theta_1$ is to define another nuclear operator, $H$, on $L^2$. Applying, if necessary, Theorem 2 to the four-element family $\mathcal{S} = \left\{H, S, A := (H^*H + S^*S)^{\frac{1}{2}}, \widetilde{A} := (HH^* + SS^*)^{\frac{1}{2}}\right\}$ of nuclear operators on $\mathcal{H} = L^2$, we may and do assume (with no loss of generality) that the function $\boldsymbol{T}_\lambda := \boldsymbol{H} - \lambda \boldsymbol{S}$ is a $K^0$ kernel of Mercer type for any fixed $\lambda \in \mathbb{C}$, and that the inducing kernels $\boldsymbol{A}$, $\widetilde{\boldsymbol{A}}$ of the above (positive nuclear) operators $A$, $\widetilde{A}$ are also $K^0$ kernels, satisfying therefore $\int \widetilde{\boldsymbol{A}}(s,s)\, ds =: \mathrm{tr}\widetilde{A} < +\infty$, $\int \boldsymbol{A}(s,s)\, ds =: \mathrm{tr}A < +\infty$ (see, e.g., [11]). With these "honorably acquired" assumptions on kernels, the main purpose of what follows is to adapt the above method for solving equation (1) to equation (4).

*6) Polynomial Fredholm Series:* For $n$, $p = 0, 1, 2 \ldots$ and $s_i$, $t_i \in \mathbb{R}$ $(1 \leqslant i \leqslant p)$, define a polynomial $B_n^p[\boldsymbol{T}_\lambda]$ in a parameter $\lambda$ of degree at most $p + n$ by

$$B_n^p[\boldsymbol{T}_\lambda]\begin{pmatrix} s_1 \ldots s_p \\ t_1 \ldots t_p \end{pmatrix} := \frac{1}{n!} \int \cdots \int \boldsymbol{T}_\lambda \begin{pmatrix} s_1 \ldots s_p\ \xi_1 \ldots \xi_n \\ t_1 \ldots t_p\ \xi_1 \ldots \xi_n \end{pmatrix} d\xi_1 \ldots d\xi_n \qquad (5)$$

where $\boldsymbol{T}_\lambda \binom{x_1 \cdots x_\nu}{y_1 \cdots y_\nu}$ is the $\nu$-th compound of $\boldsymbol{T}_\lambda$.

*Remark 1:* Let $\lambda \in \mathbb{C}$, $p \geqslant 1$, and $n$, be fixed. It is because every $\nu$-th iterant $\boldsymbol{T}_\lambda^\nu$,

$$\boldsymbol{T}_\lambda^\nu(s,t) := \int \cdots \int \boldsymbol{T}_\lambda(s,\xi_1) \cdots \boldsymbol{T}_\lambda(\xi_{\nu-1},t)\, d\xi_1 \dots d\xi_{\nu-1},$$

of the $K^0$ kernel $\boldsymbol{T}_\lambda$ is again a $K^0$ kernel that the function

$$\begin{pmatrix} s_1 \dots s_p \\ t_1 \dots t_p \end{pmatrix} \mapsto B_n^p[\boldsymbol{T}_\lambda] \begin{pmatrix} s_1 \dots s_p \\ t_1 \dots t_p \end{pmatrix}$$

is in $C\left(\mathbb{R}^{2p}, \mathbb{C}\right)$ and the functions

$$\begin{pmatrix} s_1 & \dots & s_m & s_{m+1} & \dots & s_p \\ t_1 & \dots & t_{m-1} & t_{m+1} & \dots & t_p \end{pmatrix} \mapsto B_n^p[\boldsymbol{T}_\lambda] \begin{pmatrix} s_1 & \dots & s_{m-1} & s_m & s_{m+1} & \dots & s_p \\ t_1 & \dots & t_{m-1} & \cdot & t_{m+1} & \dots & t_p \end{pmatrix},$$

$$\begin{pmatrix} s_1 & \dots & s_{m-1} & s_{m+1} & \dots & s_p \\ t_1 & \dots & t_m & t_{m+1} & \dots & t_p \end{pmatrix} \mapsto B_n^p[\boldsymbol{T}_\lambda] \begin{pmatrix} s_1 & \dots & s_{m-1} & \cdot & s_{m+1} & \dots & s_p \\ t_1 & \dots & t_{m-1} & t_m & t_{m+1} & \dots & t_p \end{pmatrix}$$

are in $C\left(\mathbb{R}^{2p-1}, L^2\right)$ for each $m$. Since $\boldsymbol{T}_\lambda$ and all its iterants are Hilbert-Schmidt kernels on $\mathbb{R}^2$, $B_n^p[\boldsymbol{T}_\lambda]$ can also be viewed as a $C\left(\mathbb{R}^{2p-2}, L^2\left(\mathbb{R}^2\right)\right)$ function if for its range space $L^2\left(\mathbb{R}^2\right)$ we choose as underlying space any of the $(s_i, t_j)$-spaces ($1 \leqslant i, j \leqslant p$).

*Theorem 3:* For $p = 0, 1, 2, \dots$, consider the series

$$D_p[\boldsymbol{T}_\lambda] \begin{pmatrix} s_1 \dots s_p \\ t_1 \dots t_p \end{pmatrix} := \sum_{n=0}^{\infty} B_n^p[\boldsymbol{T}_\lambda] \begin{pmatrix} s_1 \dots s_p \\ t_1 \dots t_p \end{pmatrix}, \tag{6}$$

which we call the *polynomial $p$-th Fredholm series* for $\boldsymbol{T}_\lambda$. Then (i) the following series converge for any compact subset $\mathfrak{K}$ in $\mathbb{C}$:

$$\sum_{n=0}^{\infty} \sup_{\lambda \in \mathfrak{K}} \|B_n^p\left[\boldsymbol{T}_\lambda\right]\|_{C(\mathbb{R}^{2p},\mathbb{C})} \quad (p \geqslant 0), \tag{7}$$

$$\sum_{n=0}^{\infty} \sup_{\lambda \in \mathfrak{K}} \|B_n^p\left[\boldsymbol{T}_\lambda\right]\|_{C(\mathbb{R}^{2p-1},L^2)} \quad (p \geqslant 1), \tag{8}$$

$$\sum_{n=0}^{\infty} \sup_{\lambda \in \mathfrak{K}} \|B_n^p\left[\boldsymbol{T}_\lambda\right]\|_{C(\mathbb{R}^{2p-2},L^2(\mathbb{R}^2))} \quad (p \geqslant 1), \tag{9}$$

meaning that the sum-function $D_p[\boldsymbol{T}_\lambda] \binom{s_1 \cdots s_p}{t_1 \cdots t_p}$ in (6), viewed (via Remark 1) as a function from $\mathbb{C}$ to $C\left(\mathbb{R}^{2p}, \mathbb{C}\right)$, to $C\left(\mathbb{R}^{2p-1}, L^2\right)$, or to $C\left(\mathbb{R}^{2p-2}, L^2\left(\mathbb{R}^2\right)\right)$, is an integral function of $\lambda$; (ii) for some $p$, the series (6) defines a sum-function $D_p[\boldsymbol{T}_\lambda] : \mathbb{C} \to C\left(\mathbb{R}^{2p}, \mathbb{C}\right)$ whose value is not $\theta_p$ at some point $\lambda \in \mathbb{C}$.

*Proof:* (i) For each $\lambda \in \mathbb{C}$, let $T_\lambda = W_\lambda \left(V_\lambda\right)^*$ be an arbitrary but fixed $\mathcal{M}$ factorization for $T_\lambda := H - \lambda S$, and let $\boldsymbol{F}_\lambda$, $\boldsymbol{G}_\lambda$ denote those $K^0$ kernels of Mercer type which induce the integral positive operators $F_\lambda := V_\lambda \left(V_\lambda\right)^*$, $G_\lambda := W_\lambda \left(W_\lambda\right)^*$, respectively. Let $\{u_n\}$ be an arbitrary but fixed orthonormal basis for $L^2$. Then for each fixed $\lambda \in \mathbb{C}$

$$\begin{aligned} \boldsymbol{T}_\lambda(s,t) &= \sum_n W_\lambda u_n(s) \overline{V_\lambda u_n(t)}, \\ \boldsymbol{F}_\lambda(s,t) &= \sum_n V_\lambda u_n(s) \overline{V_\lambda u_n(t)}, \\ \boldsymbol{G}_\lambda(s,t) &= \sum_n W_\lambda u_n(s) \overline{W_\lambda u_n(t)} \end{aligned} \tag{10}$$

for all $s$, $t \in \mathbb{R}$, where the series are $\mathbb{C}$-absolutely convergent in $C\left(\mathbb{R}^2, \mathbb{C}\right)$, by Theorem 1. For definiteness, assume that, for any $\lambda \in \mathbb{C}$, $F_\lambda = |T_\lambda|$ and $G_\lambda = \left|(T_\lambda)^*\right|$ (compare with Example 1). Then, given $f \in L^2$ and $\lambda \in \mathbb{C}$, the following chain of relations holds:

$$\begin{aligned} &\left\langle (F_\lambda)^2 f, f \right\rangle = \left\langle |T_\lambda|^2 f, f \right\rangle \\ &= \left\langle |H|^2 f, f \right\rangle - \overline{\lambda \left\langle Sf, Hf \right\rangle} - \lambda \left\langle Sf, Hf \right\rangle + |\lambda|^2 \left\langle |S|^2 f, f \right\rangle \\ &= \|Hf\|^2 - 2\mathrm{Re}\left(\lambda \left\langle Sf, Hf \right\rangle\right) + |\lambda|^2 \|Sf\|^2 \\ &\leqslant \|Hf\|^2 + 2|\lambda| \, \|Hf\| \, \|Sf\| + |\lambda|^2 \|Sf\|^2 \\ &\leqslant 2\left(\|Hf\|^2 + |\lambda|^2 \|Sf\|^2\right) \\ &\leqslant 2\left(1 + |\lambda|^2\right)\left(\|Hf\|^2 + \|Sf\|^2\right) \\ &= \left\langle 2(1+|\lambda|^2)(|H|^2 + |S|^2) f, f \right\rangle = \left\langle 2(1+|\lambda|^2) A^2 f, f \right\rangle. \end{aligned}$$

Hence, whenever $\mathfrak{K}$ is a compact subset of $\mathbb{C}$ and $c(\mathfrak{K}) := \sup_{\lambda \in \mathfrak{K}} \sqrt{2\left(1+|\lambda|^2\right)}$, the inequality $F_\lambda \leqslant c(\mathfrak{K}) A$ holds at all $\lambda$ in $\mathfrak{K}$. In a similar manner it may be proved that, for each $\lambda$ in $\mathfrak{K}$, $G_\lambda = \left|(T_\lambda)^*\right| \leqslant c(\mathfrak{K}) \left(|H^*|^2 + |S^*|^2\right)^{\frac{1}{2}} = c(\mathfrak{K}) \widetilde{A}$. In terms of kernels, these operator inequalities mean that

$$\boldsymbol{F}_\lambda(s,s) \leqslant c(\mathfrak{K}) \boldsymbol{A}(s,s), \quad \boldsymbol{G}_\lambda(s,s) \leqslant c(\mathfrak{K}) \widetilde{\boldsymbol{A}}(s,s) \tag{11}$$

whenever $s \in \mathbb{R}$ and $\lambda \in \mathfrak{K}$. Express the $p$-th compound of $\boldsymbol{T}_\lambda$ in terms of the functions $W_\lambda u_n$, $V_\lambda u_n$ used in (10), as follows:

$$\begin{aligned} \boldsymbol{T}_\lambda \begin{pmatrix} s_1 \dots s_p \\ t_1 \dots t_p \end{pmatrix} &= \sum_{\pi \in \mathbb{S}_p} \mathrm{sgn}(\pi) \boldsymbol{T}_\lambda(s_1, t_{\pi_1}) \cdots \boldsymbol{T}_\lambda(s_p, t_{\pi_p}) \\ &= \sum_{\pi \in \mathbb{S}_p} \left( \mathrm{sgn}(\pi) \prod_{i=1}^{p} \sum_n W_\lambda u_n(s_i) \overline{V_\lambda u_n\left(t_{\pi_i}\right)} \right) \\ &= \sum_{n_1,\dots,n_p=1}^{\infty} \left( \sum_{\pi \in \mathbb{S}_p} \mathrm{sgn}(\pi) \prod_{i=1}^{p} W_\lambda u_{n_i}(s_i) \overline{V_\lambda u_{n_i}\left(t_{\pi_i}\right)} \right) \\ &= \sum_{n_1,\dots,n_p=1}^{\infty} \prod_{i=1}^{p} W_\lambda u_{n_i}(s_i) \overline{V_\lambda u \begin{pmatrix} n_1 \dots n_p \\ t_1 \dots t_p \end{pmatrix}} \\ &= \sum_{n_1 < \cdots < n_p} \left( \sum_{\pi \in \mathbb{S}_p} \prod_{i=1}^{p} W_\lambda u_{n_i}\left(s_{\pi_i}\right) \overline{V_\lambda u \begin{pmatrix} n_1 \dots n_p \\ t_{\pi_1} \dots t_{\pi_p} \end{pmatrix}} \right) \\ &= \sum_{n_1 < \cdots < n_p} W_\lambda u \begin{pmatrix} n_1 \dots n_p \\ s_1 \dots s_p \end{pmatrix} \overline{V_\lambda u \begin{pmatrix} n_1 \dots n_p \\ t_1 \dots t_p \end{pmatrix}}, \end{aligned}$$

where the notation $f \binom{n_1 \cdots n_p}{s_1 \cdots s_p}$ means $\det[f_{n_i}(s_j)]_{i,j=1}^p$, $\mathbb{S}_p$ represents the symmetric group on $p$ symbols, and the last series converges $\mathbb{C}$-absolutely in $C\left(\mathbb{R}^{2p}, \mathbb{C}\right)$. The corresponding bilinear series representations for the $p$-th compounds of $\boldsymbol{F}_\lambda$ and $\boldsymbol{G}_\lambda$ are:

$$\boldsymbol{F}_\lambda \begin{pmatrix} s_1 \dots s_p \\ t_1 \dots t_p \end{pmatrix} = \sum_{n_1 < \cdot < n_p} V_\lambda u \begin{pmatrix} n_1 \dots n_p \\ s_1 \dots s_p \end{pmatrix} \overline{V_\lambda u \begin{pmatrix} n_1 \dots n_p \\ t_1 \dots t_p \end{pmatrix}},$$

$$\boldsymbol{G}_\lambda \begin{pmatrix} s_1 \dots s_p \\ t_1 \dots t_p \end{pmatrix} = \sum_{n_1 < \cdot < n_p} W_\lambda u \begin{pmatrix} n_1 \dots n_p \\ s_1 \dots s_p \end{pmatrix} \overline{W_\lambda u \begin{pmatrix} n_1 \dots n_p \\ t_1 \dots t_p \end{pmatrix}}.$$

Use them and Cauchy's inequality to arrive at a first estimate of the integrand compound in (5):

$$\begin{aligned}&\left|\boldsymbol{T}_\lambda\begin{pmatrix} s_1 \ldots s_p \; \xi_1 \ldots \xi_n \\ t_1 \ldots t_p \; \xi_1 \ldots \xi_n \end{pmatrix}\right|^2 \\ &\leqslant \boldsymbol{G}_\lambda\begin{pmatrix} s_1 \ldots s_p \; \xi_1 \ldots \xi_n \\ s_1 \ldots s_p \; \xi_1 \ldots \xi_n \end{pmatrix} \boldsymbol{F}_\lambda\begin{pmatrix} t_1 \ldots t_p \; \xi_1 \ldots \xi_n \\ t_1 \ldots t_p \; \xi_1 \ldots \xi_n \end{pmatrix}.\end{aligned} \tag{12}$$

For the next estimates, essential use will be made of Fisher's inequality which states that $\det C \leqslant \det C_0 \cdot \det C_1$, where $C$ is an $m \times m$ Hermitian matrix whose principal minors are nonnegative, $C_0$ is the $k \times k$ submatrix obtained from $C$ by retaining only the elements to the first $k$ rows and the first $k$ columns $(k < m)$, and $C_1$ is the $(m-k) \times (m-k)$ submatrix of $C$ obtained by deleting the first $k$ rows and first $k$ columns [7, Subsection II.4.4.2]. Apply Fisher's inequality, with $m = p + n$, $k = p$, to each factor in the right-hand side of (12):

$$\begin{aligned}&\left|\boldsymbol{T}_\lambda\begin{pmatrix} s_1 \ldots s_p \; \xi_1 \ldots \xi_n \\ t_1 \ldots t_p \; \xi_1 \ldots \xi_n \end{pmatrix}\right|^2 \\ &\leqslant \boldsymbol{G}_\lambda\begin{pmatrix} s_1 \ldots s_p \\ s_1 \ldots s_p \end{pmatrix}\boldsymbol{F}_\lambda\begin{pmatrix} t_1 \ldots t_p \\ t_1 \ldots t_p \end{pmatrix}\boldsymbol{G}_\lambda\begin{pmatrix} \xi_1 \ldots \xi_n \\ \xi_1 \ldots \xi_n \end{pmatrix}\boldsymbol{F}_\lambda\begin{pmatrix} \xi_1 \ldots \xi_n \\ \xi_1 \ldots \xi_n \end{pmatrix} \\ &\leqslant \boldsymbol{G}_\lambda\begin{pmatrix} s_1 \ldots s_p \\ s_1 \ldots s_p \end{pmatrix}\boldsymbol{F}_\lambda\begin{pmatrix} t_1 \ldots t_p \\ t_1 \ldots t_p \end{pmatrix}\frac{\left[\boldsymbol{G}_\lambda\begin{pmatrix} \xi_1 \ldots \xi_n \\ \xi_1 \ldots \xi_n \end{pmatrix}+\boldsymbol{F}_\lambda\begin{pmatrix} \xi_1 \ldots \xi_n \\ \xi_1 \ldots \xi_n \end{pmatrix}\right]^2}{4}.\end{aligned}$$

Next, make a multiple application of Fisher's inequality with $k = 1$ to the compounds in the right-hand side of the last inequality and then make use of (11) to get:

$$\begin{aligned}&\left|\boldsymbol{T}_\lambda\begin{pmatrix} s_1 \ldots s_p \; \xi_1 \ldots \xi_n \\ t_1 \ldots t_p \; \xi_1 \ldots \xi_n \end{pmatrix}\right| \\ &\leqslant \left[\prod_{i=1}^{p}\boldsymbol{G}_\lambda(s_i,s_i)\boldsymbol{F}_\lambda(t_i,t_i)\right]^{\frac{1}{2}} \frac{\prod_{i=1}^{n}\boldsymbol{G}_\lambda(\xi_i,\xi_i)+\prod_{i=1}^{n}\boldsymbol{F}_\lambda(\xi_i,\xi_i)}{2} \\ &\leqslant \left[\prod_{i=1}^{p}\widetilde{\boldsymbol{A}}(s_i,s_i)\boldsymbol{A}(t_i,t_i)\right]^{\frac{1}{2}} \frac{\prod_{i=1}^{n}\widetilde{\boldsymbol{A}}(\xi_i,\xi_i)+\prod_{i=1}^{n}\boldsymbol{A}(\xi_i,\xi_i)}{2c(\mathfrak{K})^{-(p+n)}}\end{aligned}$$

for all $\lambda \in \mathfrak{K}$. Integrate both ends with respect to $\xi_1, \ldots, \xi_n$ over $\mathbb{R}^n$ to derive the following inequalities valid for all $\lambda \in \mathfrak{K}$:

$$\begin{aligned}&\left|B_n^p[\boldsymbol{T}_\lambda]\begin{pmatrix} s_1 \ldots s_p \\ t_1 \ldots t_p \end{pmatrix}\right| \leqslant \\ &\left[\prod_{i=1}^{p}\widetilde{\boldsymbol{A}}\left(s_i,s_i\right)\boldsymbol{A}\left(t_i,t_i\right)\right]^{\frac{1}{2}} \frac{\left(\mathrm{tr}\widetilde{A}\right)^n+\left(\mathrm{tr}A\right)^n}{n!2c(\mathfrak{K})^{-(p+n)}} \\ &\leqslant \frac{M^{2p}c(\mathfrak{K})^{p+n}\left(\left(\mathrm{tr}\widetilde{A}\right)^n+\left(\mathrm{tr}A\right)^n\right)}{2n!}\end{aligned} \tag{13}$$

where $M = \max\left\{\sup_{s\in\mathbb{R}}\widetilde{\boldsymbol{A}}(s,s)^{\frac{1}{2}}, \sup_{t\in\mathbb{R}}\boldsymbol{A}(t,t)^{\frac{1}{2}}\right\}$; whence the series (7) converges. Passage to the $C\left(\mathbb{R}^{2p-1}, L^2\right)$ norm in the first inequality of (13) yields

$$\begin{aligned}&\|B_n^p\left[\boldsymbol{T}_\lambda\right]\|_{C(\mathbb{R}^{2p-1},L^2)} \\ &\leqslant \frac{M^{2p-1}c(\mathfrak{K})^{p+n}\max\left\{\mathrm{tr}\widetilde{A},\mathrm{tr}A\right\}^{\frac{1}{2}}\left(\left(\mathrm{tr}\widetilde{A}\right)^n+\left(\mathrm{tr}A\right)^n\right)}{2n!}\end{aligned}$$

for all $\lambda \in \mathfrak{K}$, regarless of the particular choice from $s_i$, $t_i$ $(i = 1, \ldots, p)$ of the underlying space variable for $L^2$ (cf. Remark 1). It follows that the series (8) converges. Again from (13), by Remark 1, it holds that

$$\begin{aligned}&\|B_n^p[\boldsymbol{T}_\lambda]\|_{C(\mathbb{R}^{2p-2},L^2(\mathbb{R}^2))} \\ &\leqslant \frac{M^{2p-2}c(\mathfrak{K})^{p+n}\left(\mathrm{tr}\widetilde{A}\cdot\mathrm{tr}A\right)^{\frac{1}{2}}\left(\left(\mathrm{tr}\widetilde{A}\right)^n+\left(\mathrm{tr}A\right)^n\right)}{2n!}\end{aligned}$$

for all $\lambda \in \mathfrak{K}$, which proves that the series (9) also converges.

(ii) For the proof of statement (ii), note that, by construction,

$$D_p\left[\boldsymbol{T}_{\lambda=0}\right]\begin{pmatrix} s_1 \ldots s_p \\ t_1 \ldots t_p \end{pmatrix} = \boldsymbol{H}\left(\begin{matrix} s_1 \ldots s_p \\ t_1 \ldots t_p \end{matrix}\middle| -1\right)$$

for all $p$ (cf. (5), (6) and (2)). Hence there is a non-negative integer $p$ such that $D_p[\boldsymbol{T}_{\lambda=0}] \neq \theta_p$ (cf. (3)). ∎

*7) Solution of Equation* (4)*:* The goal of the following theorem is to present a method for solving integral equation (4), which involves the use of the polynomial Fredholm series (6) for its kernel. Let $D_p^{(j)}\left[\boldsymbol{T}_\lambda\right]$ denote the $j$-th derivative, with respect to $\lambda$, of the sum-function $D_p\left[\boldsymbol{T}_\lambda\right] : \mathbb{C} \to C\left(\mathbb{R}^{2p}, \mathbb{C}\right)$ of the polynomial $p$-th Fredholm series (6), which (in view of Theorem 3) is a $C\left(\mathbb{R}^{2p}, \mathbb{C}\right)$-valued integral function of $\lambda$.

*Theorem 4:* Let $\lambda_0 \in \mathbb{C}$ be arbitrary but fixed. Then

(I) there are unique non-negative integers $\mathbf{d} = \mathbf{d}(\lambda_0)$, $\mathbf{r} = \mathbf{r}(\lambda_0)$ such that $D_{\mathbf{d}}^{(\mathbf{r})}[\boldsymbol{T}_{\lambda_0}] \neq \theta_{\mathbf{d}}$, $D_p^{(k)}[\boldsymbol{T}_{\lambda_0}] = \theta_p$ whenever $0 \leqslant k < \mathbf{r}$ and $0 \leqslant p$, and such that $D_p^{(\mathbf{r})}[\boldsymbol{T}_{\lambda_0}] = \theta_p$ whenever $0 \leqslant p < \mathbf{d}$;

(II) if, in addition, the point $\begin{pmatrix} s'_1 \ldots s'_{\mathbf{d}} \\ t'_1 \ldots t'_{\mathbf{d}} \end{pmatrix} \in \mathbb{R}^{2\mathbf{d}}$ satisfies $D_{\mathbf{d}}^{(\mathbf{r})}\left[\boldsymbol{T}_{\lambda_0}\right]\begin{pmatrix} s'_1 \ldots s'_{\mathbf{d}} \\ t'_1 \ldots t'_{\mathbf{d}} \end{pmatrix} =: \delta \neq \theta_0$, then the functions

$$\phi_i(s) = D_{\mathbf{d}}^{(\mathbf{r})}\left[\boldsymbol{T}_{\lambda_0}\right]\begin{pmatrix} s'_1 \ldots s'_{i-1} \; s \; s'_{i+1} \ldots s'_{\mathbf{d}} \\ t'_1 \ldots t'_{i-1} \; t'_i \; t'_{i+1} \ldots t'_{\mathbf{d}} \end{pmatrix}$$

$(s \in \mathbb{R},\ 1 \leqslant i \leqslant \mathbf{d})$ form a basis for the set of solutions $f$ of the homogeneous equation $f(s) + \int \boldsymbol{T}_{\lambda_0}(s,t)f(t)\,dt = 0$, and the functions

$$\psi_l(t) = \overline{D_{\mathbf{d}}^{(\mathbf{r})}\left[\boldsymbol{T}_{\lambda_0}\right]\begin{pmatrix} s'_1 \ldots s'_{l-1} \; s'_l \; s'_{l+1} \ldots s'_{\mathbf{d}} \\ t'_1 \ldots t'_{l-1} \; t \; t'_{l+1} \ldots t'_{\mathbf{d}} \end{pmatrix}}$$

$(t \in \mathbb{R}, 1 \leqslant l \leqslant \mathbf{d})$ form a basis for the set of solutions $f$ of the conjugate homogeneous equation $f(t) + \int \overline{\boldsymbol{T}_{\lambda_0}(s,t)}f(s)\,ds = 0$. The equation (4), $f(s) + \int \boldsymbol{T}_{\lambda_0}(s,t)f(t)\,dt = g(s)$ for a.e. $s$ in $\mathbb{R}$, is solvable if and only if $\langle g, \psi_l\rangle = 0$ $(1 \leqslant l \leqslant \mathbf{d})$, and the general solution is then given by

$$f(s) = g(s) - \int\frac{D_{\mathbf{d}+1}^{(\mathbf{r})}[\boldsymbol{T}_{\lambda_0}]\begin{pmatrix} s \; s'_1 \ldots s'_{\mathbf{d}} \\ t \; t'_1 \ldots t'_{\mathbf{d}} \end{pmatrix}}{\delta}g(t)dt + \sum_{i=1}^{\mathbf{d}} c_i\phi_i(s)$$

where $c_1, \ldots, c_{\mathbf{d}} \in \mathbb{C}$ are arbitrary constants;

(III) $\mathbf{r}(\lambda_0) = 0$, and, if $\mathbf{m} := \inf_{\lambda\in\mathbb{C}}\mathbf{d}(\lambda)$ and $D_{\mathbf{m}}\left[\boldsymbol{T}_{\lambda_0}\right] \neq \theta_{\mathbf{m}}$, then $\mathbf{d}(\lambda_0) = \mathbf{m}$.

*Proof:* (I) By Theorem 3, there are $p \geqslant 0$ and $\lambda \in \mathbb{C}$ such that $\theta_p \neq D_p[\boldsymbol{T}_\lambda] = \sum_{n=0}^{\infty}\frac{D_p^{(n)}[\boldsymbol{T}_{\lambda_0}]}{n!}(\lambda - \lambda_0)^n$. Hence there is a non-negative integer $\nu$ with $D_p^{(\nu)}[\boldsymbol{T}_{\lambda_0}] \neq \theta_p$. The smallest non-negative integer $\nu$ such that $D_p^{(\nu)}[\boldsymbol{T}_{\lambda_0}] \neq \theta_p$ for some non-negative integer $p$ is the required number $\mathbf{r} = \mathbf{r}(\lambda_0)$, and the smallest non-negative integer $p$ such that $D_p^{(\mathbf{r})}[\boldsymbol{T}_{\lambda_0}] \neq \theta_p$ is the required number $\mathbf{d} = \mathbf{d}(\lambda_0)$. (*Note*: The argument resembles that for Theorem 5.2.1 in [15].)

(II) The proof of statement (II) is rather lengthy for inclusion here. It can be read off the proof of Theorem in [12]. The basic properties of polynomial Fredholm series needed for that proof are those established in Theorem 3. (*Note*: In [12], the $K^0$ kernels $\boldsymbol{H}$, $\boldsymbol{S}$ are not restricted to be Hilbert-Schmidt, and every polynomial $p$-th Fredholm series is postulated to be convergent in $C\left(\mathbb{R}^{2p}, \mathbb{C}\right)$ and in $C\left(\mathbb{R}^{2p-1}, L^2\right)$, uniformly with respect to $\lambda$ on compact subsets of a domain $D \subseteq \mathbb{C}$.)

(III) That $\mathbf{r}(\lambda_0) = 0$ comes from the solution procedure for equation (1) described in Subsection I-4 and from the defining properties of $\mathbf{r}(\lambda_0)$ given in statement (I). Therefore, since $D_{\mathbf{m}}[\boldsymbol{T}_\lambda] \not\equiv \theta_{\mathbf{m}}$ in $\mathbb{C}$ and $D_{\mathbf{m}}[\boldsymbol{T}_{\lambda_0}] \neq \theta_{\mathbf{m}}$, the number $\mathbf{d}(\lambda_0)$ cannot exceed $\mathbf{m}$. (*Note*: Statement (III) accords with Theorem on Holomorphic Operator-Function from [3, Theorem 5.1, p. 39].) ■

## References

[1] T. Carleman, "Zur Theorie der linearen Integralgleichungen," *Math. Z.*, vol. 9, pp. 196–217, 1921.

[2] I. Fredholm, "Sur une classe d'équations fonctionnelles," *Acta Math.*, vol. 27, pp. 365–390, 1903.

[3] I. Ts. Gohberg and M. G. Krein, *Introduction to the Theory of Linear Non-Selfadjoint Operators in Hilbert Space*. Moscow: Nauka, 1965, (in Russian).

[4] P. Halmos and V. Sunder, *Bounded Integral Operators on $L^2$ Spaces*. Berlin: Springer, 1978.

[5] V. B. Korotkov, *Integral Operators*. Novosibirsk: Nauka, 1983, (in Russian).

[6] W. V. Lovitt, *Linear Integral Equations*. New York: Dover, 1950.

[7] M. Marcus and H. Minc, *A Survey of Matrix Theory and Matrix Inequalities*. New York: Dover Publications Inc., 1992.

[8] J. Mercer, "Functions of positive and negative type, and their connection with the theory of integral equations," *Philos. Trans. Roy. Soc. London Ser. A*, vol. 209, pp. 415–446, 1909.

[9] S. G. Mikhlin, "On the convergence of Fredholm series," *Doklady AN SSSR*, vol. XLII, no. 9, pp. 374–377, 1944, (in Russian).

[10] I. M. Novitskii, "Integral representations of linear operators by smooth Carleman kernels of Mercer type," *Proc. Lond. Math. Soc. (3)*, vol. 68, no. 1, pp. 161–177, 1994.

[11] ——, "Fredholm minors for completely continuous operators," *Dal'nevost. Mat. Sb.*, vol. 7, pp. 103–122, 1999, (in Russian).

[12] ——, "Fredholm formulae for kernels which are linear with respect to parameter," *Dal'nevost. Mat. Zh.*, vol. 3, no. 2, pp. 173–194, 2002, (in Russian).

[13] ——, "Integral operators with infinitely smooth bi-Carleman kernels of Mercer type," *Int. Electron. J. Pure Appl. Math.*, vol. 2, no. 1, pp. 43–73, 2010.

[14] ——, "Kernels of integral equations can be boundedly infinitely differentiable on $\mathbb{R}^2$," in *Proc. 2011 World Congress on Engineering and Technology (CET 2011)*, vol. 2, IEEE Press, 2011, pp. 789–792.

[15] A. F. Ruston, *Fredholm Theory in Banach Spaces*. Cambridge: Cambridge Univ. Press, 1986.

[16] F. Smithies, "The Fredholm theory of integral equations," *Duke Math. J.*, vol. 8, pp. 107–130, 1941.